\documentclass[12pt,a4paper]{article}

\usepackage{amsfonts,amsmath,amssymb,amscd}
\usepackage{colordvi}

\newcommand{\ZZ}{\mathbb{Z}}
\newcommand{\QQ}{\mathbb{Q}}

\newcommand{\CC}{\mathbb{C}}

\newcommand{\Sp}{{\rm Sp}}

\newcommand{\A}{\mathcal{A}}
\newcommand{\Mat}{{\rm Mat}}
\newcommand{\calM}{\mathcal{M}}

\newcommand{\mat}[4]
{\begin{pmatrix}  #1 & #2 \\  #3 & #4 \end{pmatrix}}
\newcommand{\smat}[4]{\left( \smallmatrix #1 & #2 \\ #3 & #4
    \endsmallmatrix \right)}

\newcommand{\IM}[1]{{\rm Im}\, #1}

\renewcommand{\H}{\mathcal{H}}

\newcommand{\End}{{\rm End}}

\newcommand{\Jac}{{\rm Jac}}

\newcommand{\blockcomment}[1]{}
\renewcommand{\choose}[2]{\left(\! \smallmatrix #1 \\ #2 \endsmallmatrix \! \right)}

\begin{document}

\title{Computing Humbert Surfaces}
\author{David Gruenewald %
\thanks{Partially supported by an APA scholarship at the University of Sydney.}
}

\maketitle

\abstract

We describe an algorithm which computes components of
Humbert surfaces in terms of Rosenhain
invariants, based on Runge's method \cite{Runge99}.

\section{Background}
For general properties of abelian varieties we refer the reader to \cite{BirkLange04}.
Denote by $\H_2$ the Siegel upper half plane of degree $2$, which by
definition is the set
$$
\H_2 = \{ \tau \in \Mat_{2\times 2}(\CC)\ |\ {^t\!\tau}=\tau\ ,\ \IM(\tau)>0\}\,.
$$
Each $\tau\in \H_2$ corresponds to a principally polarized complex abelian surface $A_\tau$ with period matrix $(\tau\ I_2) \in \Mat_{2\times 4}(\CC).$  
Two abelian surfaces $A_\tau$ and $\A_\tau'$ are isomorphic if and only if there is a symplectic matrix $M=\smat{a}{b}{c}{d}\in\Sp_4(\ZZ)$ such that $\tau=M(\tau):=(a\tau+b)(c\tau+d)^{-1}$.  Quotienting out by this action, we obtain the moduli space $\A_2=\Sp_{4}(\ZZ)\backslash\H_2$ of isomorphism classes of principally polarized abelian surfaces.  It is a quasi-projective variety of dimension $3$ and is called the \textit{Siegel modular threefold}.

The sets of abelian surfaces having the same endomorphism ring form
subvarieties of $\A_2$.   Let $A$ be a principally polarised abelian
surface.  Then $\End(A)$ is an order in $\End(A)\otimes\QQ$ which is isomorphic to either a quartic CM field, an indefinite quaternion algebra, a real quadratic field or in the generic case $\QQ$.  The irreducible components of the corresponding moduli spaces in $\A_2$ which have ``extra endomorphisms'' have dimensions $0,1,2$ and are known as CM points, Shimura curves and Humbert surfaces respectively.  

Humbert \cite{HumbertWorks} showed that for each positive discriminant $\Delta$ there is a unique irreducible Humbert surface $H_\Delta$ in $\A_2$, and any matrix $\smat{\tau_1}{\tau_2}{\tau_2}{\tau_3} \in \H_2$ satisfying the equation
\begin{equation} \label{humbert-complex}
k\tau_1 + \ell\tau_2 - \tau_3 = 0
\end{equation}
lies on the Humbert surface $H_\Delta$ of discriminant $\Delta = 4k+\ell>0$.  

The function field of $\A_2$ is $\CC(j_1,j_2,j_3)$ where the $j_i(\tau)$ are the absolute Igusa invariants,
so there is an irreducible polynomial $H_\Delta(j_1,j_2,j_3)$ whose zero set is the Humbert surface of discriminant $\Delta$.  Unfortunately, working with Igusa invariants is impractical
due to the enormous degrees and coefficients of the polynomial.  One fares better by working in a finite cover of the moduli space, adding some level structure.  Runge \cite{Runge99} constructed an algorithm to compute Humbert components in the cover $\Gamma^*(2,4)\backslash \H_2$ using theta functions and their Fourier expansions.  The purpose of this paper is to apply Runge's method to $\A_2(2)$, the Siegel modular threefold with level-$2$ structure, using Rosenhain invariants.

\section{Level-$2$ structure}
Torelli's theorem implies that the map sending a curve $C$ to its Jacobian variety $\Jac(C)$ defines a birational map from the moduli space of genus 2 curves denoted $\calM_2$, to $\A_2$.
Given a genus 2 curve $y^2=\prod_{i=1}^6(x-u_i)$ over the complex numbers, we can send three of the $u_i$ to $0,1,\infty$ via a fractional linear transformation to get an isomorphic curve with a \textit{Rosenhain model}: $$y^2=x(x-1)(x-\lambda_1)(x-\lambda_2)(x-\lambda_3).$$  The $\lambda_i$ are called \textit{Rosenhain invariants}. 

The ordered tuple $(0,1,\infty,\lambda_1,\lambda_2,\lambda_3)$ determines an ordering of the Weierstrass points and a level 2 structure on the corresponding Jacobian, that is, determines a point of $\A_2(2)$.

Let $\calM_2(2)$ denote the moduli space of genus 2 curves together with a full level 2 structure. The points of $\calM_2(2)$ are given by triples $(\lambda_1,\lambda_2,\lambda_3)$ where the $\lambda_i$ are all distinct and different from $0$ and $1$.
The forgetful morphism $\calM_2(2) \rightarrow \calM_2$ is a Galois covering of degree $720=|S_6|$ where $S_6$ acts on the Weierstrass 6-tuple by permutations, followed by renormalising the first three coordinates to $(0,1,\infty)$.

As functions on $\calM_2(2)$, the Rosenhain invariants generate the coordinate ring of $\calM_2(2)$ and hence generate the function field of $\A_2(2)$.

\section{Theta constants and Rosenhain invariants}

Let $\tau \in \H_2$ and write
 $m'=(a,b)$ and $m''=(c,d)$.  The \textit{classical theta constants} (of half integral characteristic) are defined by

\begin{eqnarray*}
  \theta_{abcd}(\tau) = \sum_{x \in \ZZ^2} \mathrm{exp}\,2 \pi i\left(\frac{1}{2}(x + \frac{m'}{2})\cdot\tau
  \cdot {}^t(x + \frac{m'}{2}) + (x + \frac{m'}{2})\cdot {}^t(\frac{m''}{2})\right)
\end{eqnarray*}
where $a,b,c,d$ are either $0$ or $1$.
Write
\begin{eqnarray*}
  \vartheta_1  &=& \theta_{0000}(\tau) \\
  \vartheta_2  &=& \theta_{0011}(\tau) \\
  \vartheta_3  &=& \theta_{0010}(\tau) \\
  \vartheta_4  &=& \theta_{0001}(\tau) \\
  \vartheta_8  &=& \theta_{1100}(\tau) \\
  \vartheta_{10} &=& \theta_{1111}(\tau) \,.
\end{eqnarray*}
These are the only theta constants we shall be using.  
As a function of $\tau\in\A_2$ there are $720$ different Rosenhain invariant triples, any of which may be used.  We use the same choice of Rosenhain triple that Gaudry uses in \cite{Gaudry05}:
\begin{equation*}
  e_1 = \frac{\vartheta_1^2 \vartheta_3^2}{\vartheta_2^2
  \vartheta_4^2},\ 
  e_2 = \frac{\vartheta_3^2 \vartheta_8^2}{\vartheta_4^2
  \vartheta_{10}^2},\ 
  e_3 = \frac{\vartheta_1^2 \vartheta_8^2}{\vartheta_2^2
  \vartheta_{10}^2}\,.
\end{equation*}

\section{Fourier series expansions}

Let us describe the Fourier expansion of even theta constants restricted
to a Humbert surface of discriminant $\Delta \equiv 0$ or  $1$ mod
$4$, adapted from ideas in Runge's paper \cite{Runge99}. Write
$\Delta=4k + \ell$ where $\ell$ is either $0$ or $1$, and $k$ is uniquely determined.
From equation (\ref{humbert-complex}) the Humbert surface of discriminant $\Delta$ can be defined by the set
\begin{equation*}
  H_\Delta = \left\{\mat{\tau_1}{\tau_2}{\tau_2}{k\tau_1 + \ell\tau_2}
  \in \H_2 \right\}
\end{equation*}
modulo the usual $\Sp_4(\ZZ)$ equivalence relation.  Restrict
$\theta_{abcd}$ to $H_\Delta$ to get
\begin{eqnarray*}
    \theta_{abcd}(\tau) 
&=& \sum_{(x_1, x_2) \in \ZZ^2} e^{\pi i(x_1c+x_2d)}  r^{(2x_1+a)^2 + k(2x_2+b)^2}  q^{2(2x_1+a)(2x_2+b) + \ell(2x_2+b)^2}
\end{eqnarray*}
where $r=e^{2\pi i\tau_1/8}$ and $q=e^{2\pi i\tau_2/8}$.
Unfortunately, $q$ has negative exponents which computationally makes it difficult to
work with this expansion.  To overcome this difficulty, make the
invertible substitution $r=pq$ to produce the expansion
\begin{eqnarray*}
  \sum_{(x_1, x_2) \in \ZZ^2} (-1)^{x_1c+x_2d}  p^{(2x_1+a)^2 +
  k(2x_2+b)^2}  q^{(2x_1+a + 2x_2+b)^2 + (k+\ell-1)(2x_2+b)^2}
\end{eqnarray*}
which is more computationally friendly, being a power series with
integer coefficients.  Call the above expansion the \textit{Fourier expansion of
$\theta_{abcd}$ restricted to $H_\Delta$}.

Addition and muliptication of restricted Fourier expansions are just
the usual addition and muliptication operations in $\ZZ[[p,q]]$.  To
compute the expansions of Rosenhain invariants we need to know how to
invert elements of $\ZZ[[p,q]]$ where possible.  It is well known fact
about power series rings that
if $f(p,q)$ is in $\ZZ[[p,q]]$ with $f(0,0)\neq 0$ , then $f(p,q)$ is a
unit with inverse given by the geometric series
\begin{equation*}
  f(0,0)^{-1}\sum_{n\geq 0}\left(1-\frac{f(p,q)}{f(0,0)}\right)^n.
\end{equation*}
An implementation on a computer uses truncated Fourier
expansions, where arithmetic is done in $\ZZ[[p,q]]/(p^N,q^N)$ for
some positive $N$.  It is easy to see that the geometric ratio has
zero constant term, in particular $(1-f/f(0,0))^k \in (p^N,q^N)$ for
$k \geq N$ so the above formula converges to the truncated expansion
of $f^{-1}$ for any chosen precision.

From the expansions we observe that
$\vartheta_1,\vartheta_2,\vartheta_3,\vartheta_4$ have constant term
$1$, hence are invertible, but
$\vartheta_8 = 2p^{1+k}q^{k+\ell-1} + \ldots$ and $\vartheta_{10} =
-2p^{1+k}q^{k+\ell-1} + \ldots$ have zero constant term.  Fortunately one can
show that $\vartheta_8,\vartheta_{10}$ are in the ideal
$(p^{1+k}q^{k+\ell-1})\ZZ[[p,q]]$ hence by cancelling out the
$p^{1+k}q^{k+\ell-1}$ factors, the quotient $\vartheta_8 / \vartheta_{10}$ makes
sense in $\ZZ[[p,q]]$.  Thus we are able to compute the Rosenhain invariants
$\lambda_1,\lambda_2,\lambda_3$ as
Fourier expansions restricted to a Humbert surface.

\section{The algorithm}

We describe an algorithm to find the equation of all irreducible
components of $\H_\Delta$ representable in terms of Rosenhain invariants.

Much arithmetic-geometric information is known about Humbert surfaces, and more generally
Hilbert modular surfaces (see \cite{Hirz-vdG},\cite{vdG}).  All Humbert
components are hypersurfaces in $\mathcal{A}_2(2)$ of the same degree.  Of use for
us is the degree $m(\Delta)$ of the defining (irreducible) polynomial
$F_{\Delta,i}$ of a Humbert component as well as the number of
components of $H_\Delta$. 
The number of Humbert components in the Satake compactification $\mathcal{A}^*_2(2)$ of $\mathcal{A}_2(2)$ is (see \cite{Besser98})
\begin{equation*}
  m(\Delta) = 
    \begin{cases}
      10 &   \textrm{ if } \Delta \equiv 1 \textrm{ mod } 8 \\
      15 &   \textrm{ if } \Delta \equiv 0 \textrm{ mod } 4 \\
      6  &   \textrm{ if } \Delta \equiv 5 \textrm{ mod } 8 \\
    \end{cases}.
\end{equation*}
The degree of an irreducible polynomial $F^*_{\Delta,i}$
 defining a Humbert component in $\mathcal{A}^*_2(2)$ 
is given by a recursive formula 
\footnote{By working with the polynomial degree rather than the component degree, we avoid the annoyance of $H_1$ having multiplicity $2$ which would otherwise complicate the formula.}
\begin{equation*}
  a_\Delta = \sum_{x>0} m(\Delta/x^2) \deg(F^*_{(\Delta/x^2),i})
\end{equation*}
where $a_\Delta$ is the coefficient of a certain modular form of weight
$5/2$ for the group $\Gamma_0(4)$, which fortunately has a more
elementary description due to a formula of Siegel,
\begin{equation*}
  a_\Delta -24\sum_{x \in \ZZ}
  \sigma_1\left(\frac{\Delta-x^2}{4}\right) =
  \begin{cases}
    12\Delta-2 & \textrm{ if $\Delta$ is a square} \\
    0 & \textrm{ otherwise}
  \end{cases}
\end{equation*}
This provides an upper bound on the degree of the polynomials $F_{\Delta,i}(e_1,e_2,e_3)$.  From computational evidence it appears $\deg F_\Delta = \deg F^*_\Delta$ for nonsquare discriminants $\Delta$ and that $\deg F_{n^2} = (1-\frac{1}{n})\deg F^*_{n^2}$ for all $n$.

\begin{table}[htbp]
  \centering
  \begin{tabular}{c|c|c|c|c|c|c|c|c|c|c|c|c}
    $\Delta$&1&4&5&8&9&12&13&16&17&20&21&24\\
    \hline
    $\deg(F^*_{\Delta,i})$&1&4&8&8&24&16&40&32&48&32&80&48\\
  \end{tabular}
  \caption{Table of degrees}
\normalsize
\end{table}

The algorithm is very simple. We have $e_1,e_2,e_3$
represented as truncated power series.  We know the degree of the
relation we are searching for.  To find an algebraic relation 
of degree $d$, compute all monomials in $e_1,e_2,e_3$ of degree atmost
$d$ and use linear algebra to find linear dependencies between
the monomials.
Once one component has been determined, the others can easily be found by
looking at the Rosenhain $S_6$-orbit of a component.

The fixed groups of the Humbert components in this model can be computed.  
As we know, $S_6$ acts on the Rosenhain invariants via the natural
action on $(0,1,\infty,e_1,e_2,e_3)$.  Let $h_\Delta$ be
the Humbert component computed using the above algorithm.  The fixed group of
$h_\Delta$ for even discriminant splits into two cases,

\begin{equation*}
  {\rm Fix}_{S_6}(h_{4k}) = 
  \begin{cases}
      G &   \textrm{ if $k$ is odd} \\
      g^{-1}Gg &   \textrm{ if $k$ is even} \\
  \end{cases}
\end{equation*}
where $G \subset S_6$ is a group of order $48$ generated by three elements
\begin{equation*}
  (0, e_1, e_3, \infty, e_2, 1),\ 
  (e_1, e_2) \textrm{ and } 
  (1, e_1, e_3, e_2);
\end{equation*}
the conjugating element is $g=(1, \infty)(e_1, e_2, e_3)$.  
Ignoring discriminant $1$ which is a special case, the fixed group of
$\Delta \equiv 1 \pmod{8}$ is a group of order $72$ generated by
\begin{equation*}
  (0, e_1)(1, e_2)(\infty, e_3),\ 
  (1, \infty),\ 
  (e_1, e_2) \textrm{ and }
  (e_2, e_3).
\end{equation*}
For $\Delta \equiv 5 \pmod{8}$ the fixed group is a group of order
$120$ generated by
\begin{equation*}
  (0, e_1)(1, e_2)(\infty, e_3),\ 
  (1, e_3, e_2, e_1, \infty) \textrm{ and }
  (\infty, e_1, e_3, e_2).
\end{equation*}
By making use of some of the simpler fixed group symmetries, we can reduce the size of the linear algebra computation.  For example, the discriminant $12$ component $h_{12}$ satisfies $h_{12}(e_2,e_1,e_3)=h_{12}(e_1,e_2,e_3)$ which means we only need roughly half the number of evaluated power series since  $e_1^ae_2^be_3^c$ and $e_1^be_2^ae_3^c$ have the same coefficient.

The runtime of the algorithm is greatly affected by the $O(\choose{d+3}{3}) = O(d^3)$
monomials that need to be evaluated.  The linear algebra solution
requires finding the kernel of a matrix with $O(\choose{d+3}{3})$ rows
and in the order of $(N/4)^2$ columns where $N$ is the precision of the power
series, which gives a runtime cost is $O(d^6N^2)$.  To have any chance of
finding a unique relation, the number of monomials must be less than the
precision used, so that the runtime is at least of order $O(d^9)$.

From the table it is evident that
the degree increases with the discriminant, so as it stands this
algorithm can only find equations with small degrees.  
Besides discriminant $21$, we managed to produce Humbert components for all the discriminants listed in the above table.  This extends the equations found in the literature (\cite{HumbertWorks},\cite{Hashimoto-Murabayashi}) which go up to discriminant $8$.  See the appendix for the equation of a discriminant $12$ Humbert component we found.

\bibliographystyle{abbrv}
\bibliography{references}

\begin{thebibliography}{1}

\bibitem{Besser98}
A.~Besser.
\newblock Elliptic fibrations of {$K3$} surfaces and {QM} {K}ummer surfaces.
\newblock {\em Math. Z.}, 228(2):283--308, 1998.

\bibitem{BirkLange04}
C.~Birkenhake and H.~Lange.
\newblock {\em Complex abelian varieties}, volume 302 of {\em Grundlehren der
  Mathematischen Wissenschaften [Fundamental Principles of Mathematical
  Sciences]}.
\newblock Springer-Verlag, Berlin, second edition, 2004.

\bibitem{Gaudry05}
P.~Gaudry.
\newblock Fast genus 2 arithmetic based on theta functions.
\newblock Preprint, 2005.

\bibitem{Hashimoto-Murabayashi}
K.-i. Hashimoto and N.~Murabayashi.
\newblock Shimura curves as intersections of {H}umbert surfaces and defining
  equations of {QM}-curves of genus two.
\newblock {\em Tohoku Math. J. (2)}, 47(2):271--296, 1995.

\bibitem{Hirz-vdG}
F.~Hirzebruch and G.~van~der Geer.
\newblock {\em Lectures on {H}ilbert modular surfaces}, volume~77 of {\em
  S\'eminaire de Math\'ematiques Sup\'erieures [Seminar on Higher
  Mathematics]}.
\newblock Presses de l'Universit\'e de Montr\'eal, Montreal, Que., 1981.
\newblock Based on notes taken by W. Hausmann and F. J. Koll.

\bibitem{HumbertWorks}
G.~Humbert.
\newblock Sur les fonctions ab\'eliennes singuli\`eres.
\newblock {\em \OE uvres}, II:297--401, 1936.

\bibitem{Runge99}
B.~Runge.
\newblock Endomorphism rings of abelian surfaces and projective models of their
  moduli spaces.
\newblock {\em Tohoku Math. J. (2)}, 51(3):283--303, 1999.

\bibitem{vdG}
G.~van~der Geer.
\newblock {\em Hilbert modular surfaces}, volume~16 of {\em Ergebnisse der
  Mathematik und ihrer Grenzgebiete (3) [Results in Mathematics and Related
  Areas (3)]}.
\newblock Springer-Verlag, Berlin, 1988.

\end{thebibliography}

\appendix
\section*{Appendix: Equation for discriminant 12}

$$
\begin{array}{c c p{11 cm}}
0 &=&
$
e_{2}^{4}e_{3}^{4} - 4 e_{2}^{4}e_{3}^{5} + 6 e_{2}^{4}e_{3}^{6} - 4
e_{2}^{4}e_{3}^{7} + e_{2}^{4}e_{3}^{8} - 4 e_{1}e_{2}^{3}e_{3}^{4} - 16
e_{1}e_{2}^{3}e_{3}^{5} + 40 e_{1}e_{2}^{3}e_{3}^{6} - 16 e_{1}e_{2}^{3}e_{3}^{7} - 4
e_{1}e_{2}^{3}e_{3}^{8} + 160 e_{1}e_{2}^{4}e_{3}^{4} - 160 e_{1}e_{2}^{4}e_{3}^{5} -
160 e_{1}e_{2}^{4}e_{3}^{6} + 160 e_{1}e_{2}^{4}e_{3}^{7} - 132
e_{1}e_{2}^{5}e_{3}^{3} - 272 e_{1}e_{2}^{5}e_{3}^{4} + 808 e_{1}e_{2}^{5}e_{3}^{5} -
272 e_{1}e_{2}^{5}e_{3}^{6} - 132 e_{1}e_{2}^{5}e_{3}^{7} + 384
e_{1}e_{2}^{6}e_{3}^{3} - 384 e_{1}e_{2}^{6}e_{3}^{4} - 384 e_{1}e_{2}^{6}e_{3}^{5} +
384 e_{1}e_{2}^{6}e_{3}^{6} - 256 e_{1}e_{2}^{7}e_{3}^{3} + 512
e_{1}e_{2}^{7}e_{3}^{4} - 256 e_{1}e_{2}^{7}e_{3}^{5} + 6 e_{1}^{2}e_{2}^{2}e_{3}^{4}
+ 40 e_{1}^{2}e_{2}^{2}e_{3}^{5} + 164 e_{1}^{2}e_{2}^{2}e_{3}^{6} + 40
e_{1}^{2}e_{2}^{2}e_{3}^{7} + 6 e_{1}^{2}e_{2}^{2}e_{3}^{8} - 160
e_{1}^{2}e_{2}^{3}e_{3}^{4} - 352 e_{1}^{2}e_{2}^{3}e_{3}^{5} - 352
e_{1}^{2}e_{2}^{3}e_{3}^{6} - 160 e_{1}^{2}e_{2}^{3}e_{3}^{7} - 272
e_{1}^{2}e_{2}^{4}e_{3}^{3} + 1344 e_{1}^{2}e_{2}^{4}e_{3}^{4} - 608
e_{1}^{2}e_{2}^{4}e_{3}^{5} + 1344 e_{1}^{2}e_{2}^{4}e_{3}^{6} - 272
e_{1}^{2}e_{2}^{4}e_{3}^{7} + 384 e_{1}^{2}e_{2}^{5}e_{3}^{2} - 416
e_{1}^{2}e_{2}^{5}e_{3}^{3} - 480 e_{1}^{2}e_{2}^{5}e_{3}^{4} - 480
e_{1}^{2}e_{2}^{5}e_{3}^{5} - 416 e_{1}^{2}e_{2}^{5}e_{3}^{6} + 384
e_{1}^{2}e_{2}^{5}e_{3}^{7} - 762 e_{1}^{2}e_{2}^{6}e_{3}^{2} + 1064
e_{1}^{2}e_{2}^{6}e_{3}^{3} - 348 e_{1}^{2}e_{2}^{6}e_{3}^{4} + 1064
e_{1}^{2}e_{2}^{6}e_{3}^{5} - 762 e_{1}^{2}e_{2}^{6}e_{3}^{6} + 384
e_{1}^{2}e_{2}^{7}e_{3}^{2} - 384 e_{1}^{2}e_{2}^{7}e_{3}^{3} - 384
e_{1}^{2}e_{2}^{7}e_{3}^{4} + 384 e_{1}^{2}e_{2}^{7}e_{3}^{5} - 4
e_{1}^{3}e_{2}e_{3}^{4} - 16 e_{1}^{3}e_{2}e_{3}^{5} + 40 e_{1}^{3}e_{2}e_{3}^{6} -
16 e_{1}^{3}e_{2}e_{3}^{7} - 4 e_{1}^{3}e_{2}e_{3}^{8} - 160
e_{1}^{3}e_{2}^{2}e_{3}^{4} - 352 e_{1}^{3}e_{2}^{2}e_{3}^{5} - 352
e_{1}^{3}e_{2}^{2}e_{3}^{6} - 160 e_{1}^{3}e_{2}^{2}e_{3}^{7} + 808
e_{1}^{3}e_{2}^{3}e_{3}^{3} - 608 e_{1}^{3}e_{2}^{3}e_{3}^{4} + 3696
e_{1}^{3}e_{2}^{3}e_{3}^{5} - 608 e_{1}^{3}e_{2}^{3}e_{3}^{6} + 808
e_{1}^{3}e_{2}^{3}e_{3}^{7} - 384 e_{1}^{3}e_{2}^{4}e_{3}^{2} - 480
e_{1}^{3}e_{2}^{4}e_{3}^{3} - 2208 e_{1}^{3}e_{2}^{4}e_{3}^{4} - 2208
e_{1}^{3}e_{2}^{4}e_{3}^{5} - 480 e_{1}^{3}e_{2}^{4}e_{3}^{6} - 384
e_{1}^{3}e_{2}^{4}e_{3}^{7} - 256 e_{1}^{3}e_{2}^{5}e_{3} + 1064
e_{1}^{3}e_{2}^{5}e_{3}^{2} - 608 e_{1}^{3}e_{2}^{5}e_{3}^{3} + 3696
e_{1}^{3}e_{2}^{5}e_{3}^{4} - 608 e_{1}^{3}e_{2}^{5}e_{3}^{5} + 1064
e_{1}^{3}e_{2}^{5}e_{3}^{6} - 256 e_{1}^{3}e_{2}^{5}e_{3}^{7} + 384
e_{1}^{3}e_{2}^{6}e_{3} - 416 e_{1}^{3}e_{2}^{6}e_{3}^{2} - 480
e_{1}^{3}e_{2}^{6}e_{3}^{3} - 480 e_{1}^{3}e_{2}^{6}e_{3}^{4} - 416
e_{1}^{3}e_{2}^{6}e_{3}^{5} + 384 e_{1}^{3}e_{2}^{6}e_{3}^{6} - 132
e_{1}^{3}e_{2}^{7}e_{3} - 272 e_{1}^{3}e_{2}^{7}e_{3}^{2} + 808
e_{1}^{3}e_{2}^{7}e_{3}^{3} - 272 e_{1}^{3}e_{2}^{7}e_{3}^{4} - 132
e_{1}^{3}e_{2}^{7}e_{3}^{5} + e_{1}^{4}e_{3}^{4} - 4 e_{1}^{4}e_{3}^{5} + 6
e_{1}^{4}e_{3}^{6} - 4 e_{1}^{4}e_{3}^{7} + e_{1}^{4}e_{3}^{8} + 160
e_{1}^{4}e_{2}e_{3}^{4} - 160 e_{1}^{4}e_{2}e_{3}^{5} - 160 e_{1}^{4}e_{2}e_{3}^{6} +
160 e_{1}^{4}e_{2}e_{3}^{7} - 272 e_{1}^{4}e_{2}^{2}e_{3}^{3} + 1344
e_{1}^{4}e_{2}^{2}e_{3}^{4} - 608 e_{1}^{4}e_{2}^{2}e_{3}^{5} + 1344
e_{1}^{4}e_{2}^{2}e_{3}^{6} - 272 e_{1}^{4}e_{2}^{2}e_{3}^{7} - 384
e_{1}^{4}e_{2}^{3}e_{3}^{2} - 480 e_{1}^{4}e_{2}^{3}e_{3}^{3} - 2208
e_{1}^{4}e_{2}^{3}e_{3}^{4} - 2208 e_{1}^{4}e_{2}^{3}e_{3}^{5} - 480
e_{1}^{4}e_{2}^{3}e_{3}^{6} - 384 e_{1}^{4}e_{2}^{3}e_{3}^{7} + 512
e_{1}^{4}e_{2}^{4}e_{3} - 348 e_{1}^{4}e_{2}^{4}e_{3}^{2} + 3696
e_{1}^{4}e_{2}^{4}e_{3}^{3} + 1496 e_{1}^{4}e_{2}^{4}e_{3}^{4} + 3696
e_{1}^{4}e_{2}^{4}e_{3}^{5} - 348 e_{1}^{4}e_{2}^{4}e_{3}^{6} + 512
e_{1}^{4}e_{2}^{4}e_{3}^{7} - 384 e_{1}^{4}e_{2}^{5}e_{3} - 480
e_{1}^{4}e_{2}^{5}e_{3}^{2} - 2208 e_{1}^{4}e_{2}^{5}e_{3}^{3} - 2208
e_{1}^{4}e_{2}^{5}e_{3}^{4} - 480 e_{1}^{4}e_{2}^{5}e_{3}^{5} - 384
e_{1}^{4}e_{2}^{5}e_{3}^{6} - 272 e_{1}^{4}e_{2}^{6}e_{3} + 1344
e_{1}^{4}e_{2}^{6}e_{3}^{2} - 608 e_{1}^{4}e_{2}^{6}e_{3}^{3} + 1344
e_{1}^{4}e_{2}^{6}e_{3}^{4} - 272 e_{1}^{4}e_{2}^{6}e_{3}^{5} + 160
e_{1}^{4}e_{2}^{7}e_{3} - 160 e_{1}^{4}e_{2}^{7}e_{3}^{2} - 160
e_{1}^{4}e_{2}^{7}e_{3}^{3} + 160 e_{1}^{4}e_{2}^{7}e_{3}^{4} + e_{1}^{4}e_{2}^{8} -
4 e_{1}^{4}e_{2}^{8}e_{3} + 6 e_{1}^{4}e_{2}^{8}e_{3}^{2} - 4
e_{1}^{4}e_{2}^{8}e_{3}^{3} + e_{1}^{4}e_{2}^{8}e_{3}^{4} - 132
e_{1}^{5}e_{2}e_{3}^{3} - 272 e_{1}^{5}e_{2}e_{3}^{4} + 808 e_{1}^{5}e_{2}e_{3}^{5} -
272 e_{1}^{5}e_{2}e_{3}^{6} - 132 e_{1}^{5}e_{2}e_{3}^{7} + 384
e_{1}^{5}e_{2}^{2}e_{3}^{2} - 416 e_{1}^{5}e_{2}^{2}e_{3}^{3} - 480
e_{1}^{5}e_{2}^{2}e_{3}^{4} - 480 e_{1}^{5}e_{2}^{2}e_{3}^{5} - 416
e_{1}^{5}e_{2}^{2}e_{3}^{6} + 384 e_{1}^{5}e_{2}^{2}e_{3}^{7} - 256
e_{1}^{5}e_{2}^{3}e_{3} + 1064 e_{1}^{5}e_{2}^{3}e_{3}^{2} - 608
e_{1}^{5}e_{2}^{3}e_{3}^{3} + 3696 e_{1}^{5}e_{2}^{3}e_{3}^{4} - 608
e_{1}^{5}e_{2}^{3}e_{3}^{5} + 1064 e_{1}^{5}e_{2}^{3}e_{3}^{6} - 256
e_{1}^{5}e_{2}^{3}e_{3}^{7} - 384 e_{1}^{5}e_{2}^{4}e_{3} - 480
e_{1}^{5}e_{2}^{4}e_{3}^{2} - 2208 e_{1}^{5}e_{2}^{4}e_{3}^{3} - 2208
e_{1}^{5}e_{2}^{4}e_{3}^{4} - 480 e_{1}^{5}e_{2}^{4}e_{3}^{5} - 384
e_{1}^{5}e_{2}^{4}e_{3}^{6} + 808 e_{1}^{5}e_{2}^{5}e_{3} - 608
e_{1}^{5}e_{2}^{5}e_{3}^{2} + 3696 e_{1}^{5}e_{2}^{5}e_{3}^{3} - 608
e_{1}^{5}e_{2}^{5}e_{3}^{4} + 808 e_{1}^{5}e_{2}^{5}e_{3}^{5} - 160
e_{1}^{5}e_{2}^{6}e_{3} - 352 e_{1}^{5}e_{2}^{6}e_{3}^{2} - 352
e_{1}^{5}e_{2}^{6}e_{3}^{3} - 160 e_{1}^{5}e_{2}^{6}e_{3}^{4} - 4 e_{1}^{5}e_{2}^{7}
- 16 e_{1}^{5}e_{2}^{7}e_{3} + 40 e_{1}^{5}e_{2}^{7}e_{3}^{2} - 16
e_{1}^{5}e_{2}^{7}e_{3}^{3} - 4 e_{1}^{5}e_{2}^{7}e_{3}^{4} + 384
e_{1}^{6}e_{2}e_{3}^{3} - 384 e_{1}^{6}e_{2}e_{3}^{4} - 384 e_{1}^{6}e_{2}e_{3}^{5} +
384 e_{1}^{6}e_{2}e_{3}^{6} - 762 e_{1}^{6}e_{2}^{2}e_{3}^{2} + 1064
e_{1}^{6}e_{2}^{2}e_{3}^{3} - 348 e_{1}^{6}e_{2}^{2}e_{3}^{4} + 1064
e_{1}^{6}e_{2}^{2}e_{3}^{5} - 762 e_{1}^{6}e_{2}^{2}e_{3}^{6} + 384
e_{1}^{6}e_{2}^{3}e_{3} - 416 e_{1}^{6}e_{2}^{3}e_{3}^{2} - 480
e_{1}^{6}e_{2}^{3}e_{3}^{3} - 480 e_{1}^{6}e_{2}^{3}e_{3}^{4} - 416
e_{1}^{6}e_{2}^{3}e_{3}^{5} + 384 e_{1}^{6}e_{2}^{3}e_{3}^{6} - 272
e_{1}^{6}e_{2}^{4}e_{3} + 1344 e_{1}^{6}e_{2}^{4}e_{3}^{2} - 608
e_{1}^{6}e_{2}^{4}e_{3}^{3} + 1344 e_{1}^{6}e_{2}^{4}e_{3}^{4} - 272
e_{1}^{6}e_{2}^{4}e_{3}^{5} - 160 e_{1}^{6}e_{2}^{5}e_{3} - 352
e_{1}^{6}e_{2}^{5}e_{3}^{2} - 352 e_{1}^{6}e_{2}^{5}e_{3}^{3} - 160
e_{1}^{6}e_{2}^{5}e_{3}^{4} + 6 e_{1}^{6}e_{2}^{6} + 40 e_{1}^{6}e_{2}^{6}e_{3} + 164
e_{1}^{6}e_{2}^{6}e_{3}^{2} + 40 e_{1}^{6}e_{2}^{6}e_{3}^{3} + 6
e_{1}^{6}e_{2}^{6}e_{3}^{4} - 256 e_{1}^{7}e_{2}e_{3}^{3} + 512
e_{1}^{7}e_{2}e_{3}^{4} - 256 e_{1}^{7}e_{2}e_{3}^{5} + 384
e_{1}^{7}e_{2}^{2}e_{3}^{2} - 384 e_{1}^{7}e_{2}^{2}e_{3}^{3} - 384
e_{1}^{7}e_{2}^{2}e_{3}^{4} + 384 e_{1}^{7}e_{2}^{2}e_{3}^{5} - 132
e_{1}^{7}e_{2}^{3}e_{3} - 272 e_{1}^{7}e_{2}^{3}e_{3}^{2} + 808
e_{1}^{7}e_{2}^{3}e_{3}^{3} - 272 e_{1}^{7}e_{2}^{3}e_{3}^{4} - 132
e_{1}^{7}e_{2}^{3}e_{3}^{5} + 160 e_{1}^{7}e_{2}^{4}e_{3} - 160
e_{1}^{7}e_{2}^{4}e_{3}^{2} - 160 e_{1}^{7}e_{2}^{4}e_{3}^{3} + 160
e_{1}^{7}e_{2}^{4}e_{3}^{4} - 4 e_{1}^{7}e_{2}^{5} - 16 e_{1}^{7}e_{2}^{5}e_{3} + 40
e_{1}^{7}e_{2}^{5}e_{3}^{2} - 16 e_{1}^{7}e_{2}^{5}e_{3}^{3} - 4
e_{1}^{7}e_{2}^{5}e_{3}^{4} + e_{1}^{8}e_{2}^{4} - 4 e_{1}^{8}e_{2}^{4}e_{3} + 6
e_{1}^{8}e_{2}^{4}e_{3}^{2} - 4 e_{1}^{8}e_{2}^{4}e_{3}^{3} +
e_{1}^{8}e_{2}^{4}e_{3}^{4}
$
\end{array}
$$
\normalsize

\end{document}